\documentclass[12pt,a4paper,twoside]{article}
\usepackage[T2A]{fontenc}
\usepackage[english,russian]{babel}

\fontsize{12}{21} \selectfont

\newtheorem{theorem}{Theorem}[section]
\newtheorem{lemma}{Lemma}[section]

\begin{document}


\makeatletter
\renewcommand{\@evenhead}{\hfil
{\bf Irina I. Bodrenko}}
\renewcommand{\@oddhead}{\hfil
\small{\underline{\bf\qquad\quad
 Some properties of Kaehler
submanifolds with recurrent tensor fields.\qquad}} }

\begin{center}
{
\large\bf
Some properties of Kaehler submanifolds\\
with recurrent tensor fields }
\end{center}

\begin{center}
{\bf Irina~I.~Bodrenko} \footnote{\copyright  Irina~I.~Bodrenko,
associate professor,
Department of Mathematics,\\
Volgograd State University,
University Prospekt 100, Volgograd, 400062, RUSSIA.\\
E.-mail: bodrenko@mail.ru \qquad \qquad http://www.bodrenko.com
\qquad \qquad http://www.bodrenko.org}
\end{center}

\begin{center}
{\bf Abstract}\\
\end{center}
\begin{center}
The properties of Kaehler submanifolds with recurrent the second
fundamental form in spaces of constant holomorphic sectional
curvature are being studied in this article.
\end{center}

\section*{Introduction}

Let $M^{2m+2l}$ be a Kaehler manifold of complex dimension $m+l$
$(m\geq 1, l\geq 1)$ with almost complex structure $J$ and a
Riemannian metric $\widetilde g$, $\widetilde\nabla$ be the
Riemannian connection coordinated with $\widetilde g$, $\widetilde
R$ be the Riemannian curvature tensor of manifold $M^{2m+2l}$. Let
$F^{2m}$ be a Kaehler submanifold of complex dimension $m$ in
$M^{2m+2l}$ with induced Riemannian metric $g$. The restriction
$J$ to $F^{2m}$ defines induced almost complex structure on
$F^{2m},$ which we will denote by the same symbol $J$. Let
$\nabla$ be the Riemannian connection coordinated with $g$, $D$ be
 the normal connection, $b$ be the second fundamental form, $R^\bot$ be
 the tensor of normal curvature
  of submanifold $F^{2m}$, $\overline\nabla=
\nabla\oplus D$ be the connection of van der Waerden ---
Bortolotti. $b$ is called {\it parallel} if $\overline\nabla
b\equiv 0$. A tensor of normal curvature $R^\bot$ is called {\it
parallel} if $\overline\nabla R^\bot \equiv 0$.

According to the definition of recurrent tensor field (see
 [1] , note 8), nonzero form $b\ne 0$ is called {\it
recurrent} if there exists 1-form $\mu$ on $F^{2m}$ such that
$\overline\nabla b = \mu\otimes b$.
\bigskip

\begin{theorem}
\label{theorem 1}. Let $F^{2m}$ be a Kaehler submanifold of
complex dimension $m$ in a Kaehler manifold $M^{2m+2l}(c)$ of
complex dimension $m+l$ and constant holomorphic sectional
curvature $c$. If $F^{2m}$ has recurrent the second fundamental
form $b$ then the tensor of normal curvature $R^\bot\ne 0$ is
parallel.
\end{theorem}

\bigskip

It is known (see [1] , note 8, theorem 3), that for a Riemannian
manifold $M$ with recurrent tensor of Riemannian curvature
$\widetilde R$ and irreducible narrowed linear group of holonomy,
it is necessary that the tensor of Riemannian curvature
$\widetilde R$ be parallel (i.e. $\widetilde\nabla\widetilde R
\equiv 0$) with the condition $\dim M \geq 3$. A Riemannian
manifold $M$ is called {\it locally symmetric} if
$\widetilde\nabla \widetilde R\equiv 0$.

\bigskip

\begin{theorem}
\label{theorem 2}. Let $F^{2m}$ be a Kaehler submanifold of
complex dimension $m$ in a Kaehler manifold $M^{2m+2l}(c)$ of
complex dimension $m+l$ and constant holomorphic sectional
curvature $c$. If $F^{2m}$ has recurrent the second fundamental
form $b$ then $F^{2m}$ is locally symmetric \\ submanifold.
\end{theorem}

\section{Main notations and formulas.}

Let $M^{n+p}$ be $(n+p)$-dimensional $(n\geq 2, p\geq 2)$ smooth
Riemannian manifold, $\widetilde g$ be a Riemannian metric on
$M^{n+p}$, $\widetilde\nabla$ be the Riemannian connection
coordinated with $\widetilde g$, $F^n$ be $n$-dimensional smooth
submanifold in $M^{n+p}$, $g$ be the induced Riemannian metric on
$F^n$, $\nabla$ be the Riemannian connection on $F^n$ coordinated
with $g$, $TF^n$ and $T^\bot F^n$ be tangent and normal bundles on
$F^n,$ respectively, $R$ and $R_1$ be the tensors of Riemannian
and Ricci curvature of connection $\nabla,$ respectively, $b$ be
the second fundamental form $F^n$, $D$ be the normal connection,
$R^\bot$ be the tensor of normal curvature, $\overline\nabla$ be
 the connection of Van der Waerden --- Bortolotti.

The formulas of Gauss and Weingarten have, respectively, the
following
 form [2] :
$$
\widetilde\nabla_X Y = \nabla_X Y + b(X, Y),
\eqno (1.1)
$$
$$
\widetilde\nabla_X\xi = -A_\xi X + D_X\xi,
\eqno (1.2)
$$
for any vector fields $X, Y,$ tangent to $F^n$, and vector field
$\xi$ normal to $F^n$.

The equations of Gauss, Peterson --- Codacci and Ricci have,
respectively, the following form [2] :
$$
\widetilde R(X, Y, Z, W) = R (X, Y, Z, W) +
\widetilde g(b(X, Z), b(Y, W)) -\widetilde g(b(X, W), b(Y, Z)),
\eqno (1.3)
$$
$$
(\widetilde R (X, Y)Z)^\bot
=(\overline\nabla_X b)(Y, Z) - (\overline\nabla_Y b)(X, Z),
\eqno (1.4)
$$
$$
\widetilde R (X, Y, \xi,\eta)
= R^\bot(X, Y, \xi, \eta) - g([A_\xi, A_\eta]X, Y),
\eqno (1.5)
$$
for any vector fields $X, Y, Z, W$, tangent to $F^n$, and vector
fields $\xi, \eta$, normal to $F^n$.

For any vector field $\xi$ normal to $F^n$, we denote as $A_\xi$
the second fundamental tensor relatively to $\xi$. For $A_\xi$ the
relation holds
$$
\widetilde g(b(X, Y), \xi) = \widetilde g(A_\xi X, Y),
\eqno (1.6)
$$
for any vector fields $X,Y$, tangent to $F^n$.

A normal vector field $\xi$ is called {\it nondegenerate} if $\det
A_\xi \ne 0$.

Covariant derivatives $\overline\nabla b$, $(\overline\nabla
A)_\xi$ and $\overline\nabla R^\bot,$ are defined by the following equalities,\\
  respectively ( [2] ):
$$
(\overline\nabla_X b)(Y, Z) = D_X (b(Y, Z)) - b(\nabla_X Y,Z)
- b(Y, \nabla_X Z),
\eqno (1.7)
$$
$$
(\overline\nabla_X A)_\xi Y= \nabla_X (A_\xi Y) - A_\xi(\nabla_X Y)
- A_{D_X \xi} Y,
\eqno (1.8)
$$
$$
(\overline\nabla_X R^\bot)(Y, Z)\xi = D_X (R^\bot(Y, Z)\xi)
- R^\bot(\nabla_X Y,Z)\xi - R^\bot(Y, \nabla_X Z)\xi -
R^\bot (Y,Z) D_X\xi,
\eqno (1.9)
$$
for any vector fields $X, Y, Z$, tangent to $F^n$, and vector
field $\xi$ normal to $F^n$.

Let indices, in this article, take the following values: $i, j, k,
s, t = 1, \dots, n$, $\alpha, \beta,\gamma = 1, \dots, p.$ We will
use the Einstein rule.

Let $x$ be an arbitrary point $F^n$, $T_xF^n$ and $T^\bot_x F^n$
be the tangent and normal spaces $F^n$ at point $x,$ respectively,
$U(x)$ be some neighborhood of point $x$, $(u^1,\dots, u^n)$ be
local coordinates on $F^n$ in $U(x)$, $\{\partial/\partial u^i\}$
be a local basis in $TF^n$, $\{n_{\alpha |}\}$ be a field of bases
of normal vectors in  $T^\bot F^n$ in $U(x)$. We may always choose
 the basis $\{n_{\alpha |}\}$ orthonormalized and assume that
$\widetilde g (n_{\alpha |}, n_{\beta|}) = \delta_{\alpha \beta}$,
where $\delta_{\alpha \beta}$ is the Kronecker symbol. We
introduce the following designations:
$$
g_{ij} =
g
\left(
\frac{\partial}{\partial u^i},\frac{\partial}{\partial u^j}
\right),
\quad
b_{ij}
=
b
\left(
\frac{\partial}{\partial u^i},\frac{\partial}{\partial u^j}
\right)
=
b^{\alpha}_{ij} n_{\alpha |},
\quad
\Gamma_{ij, k}
=
g
\left(
\nabla_{\frac{\partial}{\partial u^i}}\frac{\partial}{\partial u^j},
\frac{\partial}{\partial u^k}
\right),
$$
$$
\Gamma_{ij}^k = g^{kt}\Gamma_{ij, t},
\quad
\nabla_ib^{\alpha}_{jk}=
\frac{\partial b^{\alpha}_{jk}}{\partial u^i}
- \Gamma_{ij}^t b^\alpha_{tk} -\Gamma_{ik}^t b^\alpha_{jt},
\quad
\Gamma^\bot_{\alpha\beta |i}=
\widetilde g\left(n_{\alpha |},
\widetilde\nabla_{\frac{\partial}{\partial u^i}} n_{\beta |}\right),
$$
$$
\Gamma^{\bot\alpha}_{\beta |i}
=\delta^{\alpha\gamma} \Gamma^\bot_{\gamma\beta |i},
\quad
\Gamma^{\bot}_{\alpha\beta |i}
=\delta_{\alpha\gamma} \Gamma^{\bot\gamma}_{\beta |i},
\quad
\Gamma^{\bot}_{\alpha\beta |i} + \Gamma^{\bot}_{\beta \alpha|i} = 0,
\quad
\overline\nabla_i b^\alpha_{jk}
=\left(
\nabla_i b^\alpha_{jk}
+
\Gamma^{\bot\alpha}_{\beta |i}
b^\beta_{jk}
\right),
$$
$$
\left(
\overline\nabla_{\frac{\partial}{\partial u^i}} b
\right)
\left(
\frac{\partial}{\partial u^j},
\frac{\partial}{\partial u^k}
\right)
=
\overline\nabla_i b^\alpha_{jk} n_{\alpha |},
\qquad
b_{\alpha | ik} = \delta_{\alpha\beta} b^\beta_{ik},
\qquad
a_{\alpha| i}^j =  b_{\alpha | ik} g^{kj},
$$
$$
\nabla_i a_{\alpha| j}^k
=
\frac{\partial a_{\alpha| j}^k}{\partial u^i}
- \Gamma_{ij}^t a_{\alpha| t}^k +\Gamma_{it}^k a_{\alpha| j}^t,
\qquad
a_{\alpha| i}^j \frac{\partial}{\partial u^j}
= A_{n_{\alpha |}}
\left(
\frac{\partial}{\partial u^i}
\right),
$$
$$
\overline\nabla_i a_{\alpha| j}^k
=
\left(
\nabla_i a_{\alpha | j}^k
-
\Gamma^{\bot\beta}_{\alpha |i}
a_{\beta | j}^k
\right),
\qquad
\left(
\overline\nabla_{\frac{\partial}{\partial u^i}} A
\right)_{n_{\alpha |}}
\left(
\frac{\partial}{\partial u^j}
\right)
=
\overline\nabla_i a_{\alpha| j}^k \frac{\partial}{\partial u^k},
$$
where $\|g^{kt}\|$ and $\|\delta^{\alpha\beta}\|$ are inverse
matrixes to $\|g_{kt}\|$ and $\|\delta_{\alpha\beta}\|,$
respectively.
\medskip

We assume that a Riemannian manifold $M^{n+p}$ is almost Hermitian
manifold with almost complex structure $J$ (see [3] , chapter 6,
section 6.1). Then $M^{n+p}$ has even dimension: $n+p = 2(m+l)$,
where a number $m+l$ is called {\it complex dimension} of
$M^{n+p}$; the Riemannian metric $\widetilde g$  is almost
Hermitian, i.e. for any vector fields $\widetilde X, \widetilde
Y$, tangent to $M^{n+p}$, the following condition holds:
$$
\widetilde g (J\widetilde X, J\widetilde Y) =
\widetilde g (\widetilde X, \widetilde Y).
\eqno (1.10)
$$

Almost Hermitian manifold $M^{n+p}$ is called {\it Kaehler
manifold} ( [3] )  if almost complex structure $J$ is parallel,
i.e. for any vector fields  $\widetilde X, \widetilde Y$, tangent
to $M^{n+p}$, the following condition holds:
$$
\widetilde\nabla_{\widetilde X} J\widetilde Y =
J \widetilde\nabla_{\widetilde X}\widetilde Y.
\eqno (1.11)
$$

A submanifold $F^n$ of a Kaehler manifold $M^{n+p}$ is called
{\it Kaehler submanifold} if for any vector field $X\in TF^n,$
vector field $JX\in TF^n$. $F^n$ is Kaehler manifold relative to
induced almost complex structure $J$ and induced almost Hermitian
metric $g$ (see [3] , chapter 6, par. 6.7). Kaehler submanifold
$F^n$ in Kaehler manifold $M^{n+p},$ has even dimension $n = 2m$
and codimension $p= 2l.$ Number $m$ is called {\it complex
dimension,} and number $l$ is called {\it complex codimension} of
Kaehler submanifold $F^n$.

We denote by $M^{2m+2l}(c),$ a Kaehler manifold of complex
dimension $m+l$ of constant holomorphic sectional curvature $c$.
The tensor of Riemannian curvature $\widetilde R$ of space
$M^{2m+2l}(c)$ complies with the formula [1] :
$$
\widetilde R (\widetilde X, \widetilde Y)\widetilde Z
=
\frac{c}{4}
\left(
\widetilde g(\widetilde Y, \widetilde Z)\widetilde X
-\widetilde g(\widetilde X, \widetilde Z)\widetilde Y
+
\widetilde g(J\widetilde Y, \widetilde Z)J\widetilde X
-\widetilde g(J\widetilde X, \widetilde Z)J\widetilde Y
+ 2\widetilde g(\widetilde X, J\widetilde Y)J\widetilde Z
\right),
\eqno (1.12)
$$
for any vector fields $\widetilde X, \widetilde Y, \widetilde Z$,
tangent to $M^{2m+2l}(c)$.

\section{The properties of covariant derivative $\overline\nabla$.}

\begin{lemma}
\label{lemma 1}. Let $F^n$ be a submanifold in a Riemannian
manifold $M^{n+p}$. Then the following equality holds:
$$
\widetilde g((\overline\nabla_Z A)_\xi X, Y)
=
\widetilde g((\overline\nabla_Z b)(X, Y), \xi)
\quad
\forall X, Y, Z \in TF^n,
\quad
\forall \xi \in T^\bot F^n.
\eqno (2.1)
$$
\end{lemma}

{\bf Proof.} We will find the expressions of the left and the
right parts of the equality (2.1), in local coordinates. We assume
$$
Z =Z^i \frac{\partial}{\partial u^i},
\quad
X = X^j \frac{\partial}{\partial u^j},
\quad
Y =Y^k \frac{\partial}{\partial u^k},
\quad
\xi = \xi^\alpha n_{\alpha |}.
\eqno (2.2)
$$
We have:
$$
\widetilde g((\overline\nabla_Z A)_\xi X, Y)
=
Z^i X^j Y^k
\xi^\alpha
g_{sk}\overline\nabla_i a_{\alpha| j}^s
=
Z^i X^j Y^k
\left(
\xi^\alpha
g_{sk}
\nabla_i a_{\alpha | j}^s
-
\xi^\alpha
g_{sk}
\Gamma^{\bot\beta}_{\alpha |i}
a_{\beta | j}^s
\right)
=
$$
$$
=
Z^i X^j Y^k
\left(
\xi^\alpha
\nabla_i (g_{sk} a_{\alpha | j}^s)
-
\xi^\alpha
\Gamma^{\bot\beta}_{\alpha |i}
g_{sk} a_{\beta | j}^s
\right)
=
Z^i X^j Y^k
\left(
\xi^\alpha
\nabla_i b_{\alpha | jk}
-
\xi^\alpha
\Gamma^{\bot\beta}_{\alpha |i}
b_{\beta | jk}
\right)
=
$$
$$
=
Z^i X^j Y^k
\left(
\xi^\alpha
\nabla_i (\delta_{\alpha\beta}b^\beta_{jk})
-
\xi^\alpha
\Gamma^{\bot\beta}_{\alpha |i}
\delta_{\beta\gamma}
b^\gamma_{jk}
\right)
=
Z^i X^j Y^k
\left(
\xi^\alpha
\delta_{\alpha\beta}
\nabla_i b^\beta_{jk}
-
\xi^\alpha
\Gamma^{\bot\beta}_{\alpha |i}
\delta_{\beta\gamma}
b^\gamma_{jk}
\right)
=
$$
$$
=
Z^i X^j Y^k
\left(
\xi^\alpha
\delta_{\alpha\beta}
\nabla_i b^\beta_{jk}
-
\xi^\alpha
\Gamma^{\bot}_{\gamma\alpha |i}
b^\gamma_{jk}
\right)
=
$$
$$
=
Z^i X^j Y^k
\left(
\xi^\alpha
\delta_{\alpha\beta}
\left(
\overline\nabla_i b^\beta_{jk}
- \Gamma^{\bot\beta}_{\gamma |i}
b^\gamma_{jk}
\right)
-
\xi^\alpha
\Gamma^{\bot}_{\gamma\alpha |i}
b^\gamma_{jk}
\right)
=
$$
$$
=
Z^i X^j Y^k
\left(
\xi^\alpha
\delta_{\alpha\beta}
\overline\nabla_i b^\beta_{jk}
-
\xi^\alpha
\delta_{\alpha\beta}
\Gamma^{\bot\beta}_{\gamma |i}
b^\gamma_{jk}
-
\xi^\alpha
\Gamma^{\bot}_{\gamma\alpha |i}
b^\gamma_{jk}
\right)
=
$$
$$
=
Z^i X^j Y^k
\left(
\xi^\alpha
\delta_{\alpha\beta}
\overline\nabla_i b^\beta_{jk}
-
\xi^\alpha
\Gamma^{\bot}_{\alpha\gamma |i}
b^\gamma_{jk}
-
\xi^\alpha
\Gamma^{\bot}_{\gamma\alpha |i}
b^\gamma_{jk}
\right)
=
Z^i X^j Y^k
\delta_{\alpha\beta}
\xi^\alpha
\overline\nabla_i b^\beta_{jk}
=
$$
$$
=
\widetilde g((\overline\nabla_Z b)(X, Y), \xi).
$$
Lemma is proved.

\begin{lemma}
\label{lemma 2}. Let $F^{2m}$ be a Kaehler submanifold in a
Kaehler manifold $M^{2m+ 2l}.$ Then for any $X\in TF^{2m}$ and for
any $\xi \in T^\bot F^{2m}$ the following equality holds:
$$
\left(
\overline\nabla_X
A
\right)_{J\xi}
=
J
\left(
\overline\nabla_X
A
\right)_\xi
\eqno (2.3)
$$
\end{lemma}

{\bf Proof.} From (1.1), because of (1.11), we obtain the
following equalities (see , for example, [3] , chapter 6, section
6.1, lemma 6. 26):
$$
\nabla_X JY = J \nabla_X Y,
\quad
J b(X, Y) = b(X, JY),
\qquad
\forall X, Y\in TF^{2m}.
\eqno (2.4)
$$
From (1.2) we have:
$$
\widetilde\nabla_X J\xi = -A_{J\xi} X + D_X J\xi,
\quad
J \widetilde\nabla_X\xi = J(-A_\xi X + D_X\xi).
$$
Hence,  because of (1.11), we obtain:
$$
-A_{J\xi} X + D_X J\xi = J(-A_\xi X + D_X\xi).
$$
Therefore,
$$
-A_{J\xi} X +  J A_\xi X  = J D_X\xi - D_X J\xi.
$$
Since $F^{2m}$ is a Kaehler submanifold, then, from here, we have
$$
A_{J\xi}X = J A_\xi X,
\quad
D_X(J\xi)= J D_X\xi,
\qquad
\forall X, Y\in TF^{2m}.
\eqno (2.5)
$$
From (1.7) we have
$$
(\overline\nabla_X A)_{J\xi} Y
= \nabla_X (A_{J\xi} Y) - A_{J\xi}(\nabla_X Y)
- A_{D_X (J\xi)} Y.
$$
Hence,  using (2.4) and (2.5), we have:
$$
(\overline\nabla_X A)_{J\xi} Y
= \nabla_X J(A_\xi Y) - J A_\xi(\nabla_X Y)
- A_{J(D_X \xi)} Y
=
$$
$$
=
J \nabla_X (A_\xi Y) - J A_\xi(\nabla_X Y)
- J A_{D_X \xi} Y
= J (\overline\nabla_X A)_\xi Y.
$$
Lemma is proved.

\bigskip

\begin{lemma}
\label{lemma 3}. Let $F^{2m}$ be a Kaehler submanifold in a
Kaehler manifold $M^{2m+2l}(c)$ of constant holomorphic sectional
curvature $c$. Then  for any  $X, Y, Z \in TF^{2m}$ and for any
$\xi \in T^\bot F^{2m}$ the following equalities hold:
$$
(\overline\nabla_{JZ} b)(X,Y)
=
J
\left(
(\overline\nabla_Z b)(X,Y)
\right),
\eqno (2.6)
$$
$$
(\overline\nabla_{JZ} A)_\xi
= - J (\overline\nabla_Z A)_\xi,
\eqno (2.7)
$$
$$
J A_\xi = - A_\xi J,
\eqno (2.8)
$$
$$
J(\overline\nabla_Z A)_\xi = -(\overline\nabla_Z A)_\xi J.
\eqno (2.9)
$$
\end{lemma}

{\bf Proof.} 1. Because of (1.12), the equation (1.4) takes the
following form:
$$
(\overline\nabla_X b)(Y, Z) = (\overline\nabla_Y b)(X, Z),
\quad
\forall X, Y, Z \in TF^{2m}.
\eqno (2.10)
$$
Using (2.10), from (1.7) we obtain:
$$
(\overline\nabla_{JZ} b)(X,Y)
=
(\overline\nabla_X b)(JZ,Y)
= D_X(b(JZ,Y)) - b(\nabla_X(JZ),Y) -b(JZ, \nabla_X Y).
$$
Hence,  using (2.4) and (2.5), we have:
$$
(\overline\nabla_{JZ} b)(X,Y)
= D_X(J(b(Z,Y))) -b(J\nabla_X Z,Y) -b(JZ, \nabla_X Y)
=
$$
$$
=
J(D_X (b(Z,Y))) - J(b(\nabla_X Z,Y)) - J(b(Z, \nabla_X Y))
=
$$
$$
=
J(D_X (b(Z,Y)) - b(\nabla_X Z,Y) - b(Z, \nabla_X Y))
= J
\left(
(\overline\nabla_Z b)(X,Y)
\right).
$$
The equality (2.6) is proved.

2. Using (2.6), from (2.1) we obtain:
$$
\widetilde g((\overline\nabla_{JZ} A)_\xi X, Y)
=
\widetilde g((\overline\nabla_{JZ} b)(X, Y), \xi)
=
\widetilde g(J((\overline\nabla_Z b)(X, Y)), \xi).
$$
Hence, because of (1.10)  and equality  $J^2 = -I$, we have:
$$
\widetilde g((\overline\nabla_{JZ} A)_\xi X, Y)
=
- \widetilde g((\overline\nabla_Z b)(X, Y), J\xi)
=
-\widetilde g((\overline\nabla_Z A)_{J\xi} X, Y)
=
-\widetilde g(J((\overline\nabla_Z A)_\xi X), Y).
$$
From here we get (2.7).

3.  From (1.6), using (2.4), we obtain:
$$
\widetilde g(J A_\xi X, Y) = -\widetilde g(A_\xi X, JY)
= -\widetilde g(b(X, JY), \xi)
= -\widetilde g(b(JX, Y), \xi)
= -\widetilde g(A_\xi JX, Y).
$$
Thus,
$$
\widetilde g(J A_\xi X, Y)
= -\widetilde g(A_\xi JX, Y)
\quad
\forall X,Y\in TF^{2m},
\quad
\forall \xi \in T^\bot F^{2m}.
$$
The derived equality is equivalent to (2.8).

4. From (1.8), using (2.4) and (2.8), for any $X,Y\in TF^{2m}$ and
for any $\xi\in T^\bot F^{2m},$ we have:
$$
J
\left(
(\overline\nabla_X A)_\xi Y
\right)
=
J
\left(
\nabla_X (A_\xi Y) - A_\xi(\nabla_X Y) - A_{D_X \xi} Y
\right)
=
$$
$$
=
\nabla_X J(A_\xi Y) + A_\xi J(\nabla_X Y) + A_{D_X \xi} (JY)
=
$$
$$
=
-\nabla_X (A_\xi JY) + A_\xi (\nabla_X JY) + A_{D_X \xi} (JY)
=
-
(\overline\nabla_X A)_\xi (JY).
$$
Thus,
$$
J
\left(
(\overline\nabla_X A)_\xi Y
\right)
=
-
(\overline\nabla_X A)_\xi (JY),
\quad
\forall X,Y\in TF^{2m},
\quad
\forall \xi \in T^\bot F^{2m}.
$$
The obtained equality is equivalent to  (2.9). Lemma is proved.

\bigskip

\begin{lemma}
\label{lemma 4}. Let $F^{2m}$ be a Kaehler submanifold in a
Kaehler manifold $M^{2m+2l}.$ \\
Then the following equality holds
$$
\overline\nabla_Z
\left(
\widetilde g(X, JY)J\xi
\right) = 0
\qquad
\forall X, Y, Z \in TF^{2m},
\quad
\forall \xi \in T^\bot F^{2m}.
\eqno (2.11)
$$
\end{lemma}

{\bf Proof.} By the definition of covariant derivative
$\overline\nabla,$ we have:
$$
\overline\nabla_Z
\left(
\widetilde g(X, JY)J\xi
\right)
=
$$
$$
=
D_Z
\left(
\widetilde g(X, JY)J\xi
\right)
-
\widetilde g(\nabla_Z X, JY)J\xi
-
\widetilde g(X,\nabla_Z(JY))J\xi
-
\widetilde g(X, JY)D_Z(J\xi).
$$
We transform the right part of the last equality, writing it in
local coordinates and using the designations (2.2):
$$
\left(
\frac{\partial
\left(
g_{kl}X^k(JY)^l(J\xi)^\tau
\right)
}{\partial u^i}
+
\Gamma^{\bot\tau}_{\sigma| i}
g_{kl}X^k(JY)^l(J\xi)^\sigma
\right)
Z^i n_\tau
-
$$
$$
-
g_{kl}
\left(
\frac{\partial X^k}{\partial u^i} + \Gamma^k_{im} X^m
\right)
(JY)^l(J\xi)^\tau
Z^i n_\tau
-
g_{lk}
\left(
\frac{\partial (JY)^l}{\partial u^i} + \Gamma^l_{im} (JY)^m
\right)
X^k(J\xi)^\tau
Z^i n_\tau
-
$$
$$
-
g_{kl} X^k (JY)^l
\left(
\frac{\partial (J\xi)^\tau}{\partial u^i}
+ \Gamma^{\bot \tau}_{\sigma| i}(J\xi)^\sigma
\right)
Z^i n_\tau
=
$$
$$
=
\left(
\frac{\partial g_{kl}}{\partial u^i} X^k (JY)^l(J\xi)^\tau
-
g_{kl}\Gamma^k_{im} X^m (JY)^l(J\xi)^\tau
-
g_{lk}\Gamma^l_{im} X^k (JY)^m(J\xi)^\tau
\right)
Z^i n_\tau
=
$$
$$
=
\left(
\frac{\partial g_{kl}}{\partial u^i}
-
g_{ml}\Gamma^m_{ik}
-
g_{mk}\Gamma^m_{il}
\right)
X^k (JY)^l
(J\xi)^\tau
Z^i n_\tau
=0.
$$
Lemma is proved.
\bigskip

\begin{lemma}
\label{lemma 5}. Let $F^{2m}$ be a Kaehler submanifold in a
Kaehler manifold $M^{2m+2l}(c)$ of constant holomorphic sectional
curvature $c$. Then the following equality holds
$$
R^\bot (X, Y)\xi =
\frac{c}{2}\widetilde g(X, JY)J\xi
+ b(X, A_\xi Y) - b(Y, A_\xi X),
$$
$$
\forall X, Y \in TF^{2m},
\quad
\forall \xi \in T^\bot F^{2m}.
\eqno (2.12)
$$
\end{lemma}

{\bf Proof.} Because of (1.12), we have:
$$
\widetilde R (X,Y,\xi, \eta)=
\widetilde g (\widetilde R (X,Y)\xi, \eta)=
\frac{c}{2}\widetilde g(X, JY)\widetilde g(J\xi,\eta).
$$
Then the equation (1.5) takes the form:
$$
R^\bot(X, Y, \xi, \eta)
= \frac{c}{2}\widetilde g(X, JY)\widetilde g(J\xi,\eta)
+ \widetilde g([A_\xi, A_\eta]X, Y).
$$
We transform the second term in the right part of the obtained
equality, using self-adjointness of operator $A_\xi$:
$$
\widetilde g([A_\xi, A_\eta]X, Y)
=
\widetilde g((A_\xi A_\eta - A_\eta A_\xi)X, Y)
=
\widetilde g(A_\xi (A_\eta X), Y)
- \widetilde g(A_\eta (A_\xi X), Y)
=
$$
$$
= \widetilde g(A_\eta X, A_\xi Y) - \widetilde g(A_\xi X, A_\eta Y)
=
\widetilde g(b(A_\xi Y, X),\eta) - \widetilde g(b(A_\xi X, Y),\eta).
$$
Then for any $\eta\in T^\bot F^{2m}$ we have:
$$
R^\bot (X, Y, \xi, \eta)
\equiv
\widetilde g(R^\bot (X, Y)\xi, \eta)
=
$$
$$
=
\widetilde g(\frac{c}{2}\widetilde g(X, JY)J\xi, \eta)
+
\widetilde g(b(A_\xi Y, X),\eta) - \widetilde g(b(A_\xi X, Y),\eta).
$$
From here we obtain the equality (2.12). Lemma is proved.

\bigskip

\begin{lemma}
\label{lemma 6}. Let $F^{2m}$ be a Kaehler submanifold in a
Kaehler manifold $M^{2m+2l}(c)$ of constant holomorphic sectional
curvature $c$. Then the following equality holds
$$
(\overline\nabla_Z R^\bot)(X,Y) \xi
=
$$
$$
=
(\overline\nabla _Z b)(X, A_\xi Y)
+
b(X, (\overline\nabla_Z A)_\xi Y)
-
(\overline\nabla _Z b)(Y, A_\xi X)
-
b(Y, (\overline\nabla_Z A)_\xi X),
$$
$$
\forall X, Y, Z \in TF^{2m},
\quad
\forall \xi \in T^\bot F^{2m}.
\eqno (2.13)
$$
\end{lemma}

{\bf Proof.} From formula (1.9), using (2.12), we obtain:
$$
(\overline\nabla_Z R^\bot)(X,Y) \xi
= D_Z
\left(
\frac{c}{2}\widetilde g(X, JY)J\xi
+ b(X, A_\xi Y) - b(Y, A_\xi X)
\right)
-
$$
$$
-
\left(
\frac{c}{2}\widetilde g(\nabla_Z X, JY)J\xi
+ b(\nabla_Z X, A_\xi Y) - b(Y, A_\xi (\nabla_Z X))
\right)
-
$$
$$
-
\left(
\frac{c}{2}\widetilde g(X, J(\nabla_Z Y))J\xi
+ b(X, A_\xi (\nabla_Z Y)) - b(\nabla_Z Y, A_\xi X)
\right)
-
$$
$$
-
\left(
\frac{c}{2}\widetilde g(X, JY))J(D_Z \xi)
+ b(X, A_{D_Z\xi} Y) - b(Y, A_{D_Z\xi} X)
\right)
=
$$
$$
=
\frac{c}{2}
\Biggl(
D_Z(\widetilde g(X, JY)J\xi)
- \widetilde g(\nabla_Z X, JY)J\xi
- \widetilde g(X, J(\nabla_Z Y))J\xi
- \widetilde g(X, JY) J(D_Z \xi)
\Biggr)
+
$$
$$
+
D_Z(b(X, A_\xi Y)) - D_Z(b(Y, A_\xi X))
- b(\nabla_Z X, A_\xi Y) + b(Y, A_\xi (\nabla_Z X))
-
$$
$$
-
b(X, A_\xi (\nabla_Z Y)) + b(\nabla_Z Y, A_\xi X)
- b(X, A_{D_Z\xi} Y) + b(Y, A_{D_Z\xi} X).
$$
Hence,  using (2.3) and (2.4), we have
$$
(\overline\nabla_Z R^\bot)(X,Y) \xi
=
$$
$$
=
\frac{c}{2} \overline\nabla_Z(\widetilde g(X, JY)J\xi)
+
D_Z(b(X, A_\xi Y)) - D_Z(b(Y, A_\xi X))
- b(\nabla_Z X, A_\xi Y)
+
$$
$$
+ b(Y, A_\xi (\nabla_Z X))
-
b(X, A_\xi (\nabla_Z Y)) + b(\nabla_Z Y, A_\xi X)
- b(X, A_{D_Z\xi} Y) + b(Y, A_{D_Z\xi} X).
$$
Therefore, because of (2.11), we obtain:
$$
(\overline\nabla_Z R^\bot)(X,Y) \xi
=
\Biggl(
D_Z(b(X, A_\xi Y))
- b(\nabla_Z X, A_\xi Y)
- b(X, A_\xi (\nabla_Z Y))
- b(X, A_{D_Z\xi} Y)
\Biggr)
-
$$
$$
-
\Biggl(
 D_Z(b(Y, A_\xi X))
- b(\nabla_Z Y, A_\xi X)
- b(Y, A_\xi (\nabla_Z X))
- b(Y, A_{D_Z\xi} X)
\Biggr).
$$
Hence,  using (1.7), we obtain:
$$
(\overline\nabla_Z R^\bot)(X,Y) \xi
=
\Biggl(
(\overline\nabla_Z b)(X, A_\xi Y)
+ b(X, \nabla_Z (A_\xi Y))
- b(X, A_\xi (\nabla_Z Y))
- b(X, A_{D_Z\xi} Y)
\Biggr)
-
$$
$$
-
\Biggl(
(\overline\nabla_Z b)(Y, A_\xi X)
+ b(Y, \nabla_Z (A_\xi X))
- b(Y, A_\xi (\nabla_Z X))
- b(Y, A_{D_Z\xi} X)
\Biggr).
$$
Now, using (1.8), we obtain:
$$
(\overline\nabla_Z R^\bot)(X,Y) \xi
=
$$
$$
=
\Biggl(
(\overline\nabla_Z b)(X, A_\xi Y)
+ b(X, (\overline\nabla_Z A)_\xi Y)
\Biggr)
-
\Biggl(
(\overline\nabla_Z b)(Y, A_\xi X)
+ b(Y, (\overline\nabla_Z A)_\xi X)
\Biggr).
$$
Lemma is proved.

\begin{lemma}
\label{lemma 7}. Let $F^{2m}$ be a Kaehler submanifold in a
Kaehler manifold $M^{2m+2l}(c)$ of constant holomorphic sectional
curvature $c$. Then the following equality holds
$$
(\overline\nabla_Z R^\bot)(X,Y, \xi, \eta)
=
\widetilde g([(\overline\nabla _Z A)_\xi, A_\eta]X, Y)
+
\widetilde g([A_\xi, (\overline\nabla_Z A)_\eta] X, Y),
$$
$$
\forall X, Y, Z \in TF^{2m},
\quad
\forall \xi,\eta \in T^\bot F^{2m}.
\eqno (2.14)
$$
\end{lemma}

{\bf Proof.} Because of (2.13), we have:
$$
(\overline\nabla_Z R^\bot)(X,Y, \xi, \eta)
\equiv
\widetilde g
\left(
(\overline\nabla_Z R^\bot)(X,Y)\xi, \eta
\right)
=
\widetilde g
\left(
(\overline\nabla_Z b)(X, A_\xi Y), \eta
\right)
-
$$
$$
-
\widetilde g
\left(
(\overline\nabla_Z b)(Y, A_\xi X), \eta
\right)
+
\widetilde g
\left(
b(X, (\overline\nabla_Z A)_\xi Y), \eta
\right)
-
\widetilde g
\left(
b(Y, (\overline\nabla_Z A)_\xi X), \eta
\right).
$$
In the derived equality, we transform  the first and the second
terms using (2.1), the third and the fourth using (1.6):
$$
(\overline\nabla_Z R^\bot)(X,Y, \xi, \eta)
=
\widetilde g
\left(
(\overline\nabla_Z A)_ \eta X, A_\xi Y
\right)
-
\widetilde g
\left(
(\overline\nabla_Z A)_ \eta Y, A_\xi X
\right)
+
$$
$$
+
\widetilde g
\left(
A_\eta X, (\overline\nabla_Z A)_\xi Y)
\right)
-
\widetilde g
\left(
A_\eta Y, (\overline\nabla_Z A)_\xi X)
\right).
$$
Hence,  because of self-adjointness of operators $A_\xi$ and
$(\overline\nabla A)_\xi,$ we obtain:
$$
(\overline\nabla_Z R^\bot)(X,Y, \xi, \eta)
=
\widetilde g
\left(
A_\xi(\overline\nabla_Z A)_ \eta X,  Y
\right)
-
\widetilde g
\left(
Y, (\overline\nabla_Z A)_ \eta A_\xi X
\right)
+
$$
$$
+
\widetilde g
\left(
(\overline\nabla_Z A)_\xi A_\eta X, Y)
\right)
-
\widetilde g
\left(
Y, A_\eta (\overline\nabla_Z A)_\xi X)
\right)
=
$$
$$
=
\widetilde g
\left(
[A_\xi, (\overline\nabla_Z A)_ \eta] X,  Y
\right)
+
\widetilde g
\left(
[(\overline\nabla_Z A)_\xi, A_\eta] X, Y)
\right)
.
$$
Lemma is proved.

\begin{lemma}
\label{lemma 8}. Let $F^{2m}$ be a Kaehler submanifold in a
Kaehler manifold $M^{2m+2l}(c)$ of constant holomorphic sectional
curvature $c$. Then the following equality holds
$$
(\overline\nabla_{JZ} R^\bot)(X,Y, \xi, \eta)
=
(\overline\nabla_Z R^\bot)(X,Y, J\xi, \eta)
-
2\widetilde g([(\overline\nabla _Z A)_{J\xi}, A_\eta]X, Y),
$$
$$
\forall X, Y, Z \in TF^{2m},
\quad
\forall \xi,\eta \in T^\bot F^{2m}.
\eqno (2.15)
$$
\end{lemma}

{\bf Proof.} From (2.14) we obtain:
$$
(\overline\nabla_{JZ} R^\bot)(X,Y, \xi, \eta)
=
\widetilde g([(\overline\nabla _{JZ} A)_\xi, A_\eta]X, Y)
+
\widetilde g([A_\xi, (\overline\nabla_{JZ} A)_\eta] X, Y).
$$
Hence,  using (2.7), we have:
$$
(\overline\nabla_{JZ} R^\bot)(X,Y, \xi, \eta)
=
\widetilde g([-J(\overline\nabla _Z A)_\xi, A_\eta]X, Y)
+
\widetilde g([A_\xi, -J(\overline\nabla_Z A)_\eta] X, Y).
$$
In the derived equality, we transform the second term using (2.8)
and (2.9):
$$
[A_\xi, J(\overline\nabla_Z A)_\eta]
=
A_\xi J (\overline\nabla_Z A)_\eta - J(\overline\nabla_Z A)_\eta A_\xi
=
$$
$$
=
-J A_\xi (\overline\nabla_Z A)_\eta + (\overline\nabla_Z A)_\eta J A_\xi
=
- [J A_\xi , (\overline\nabla_Z A)_\eta].
$$
Therefore,
$$
(\overline\nabla_{JZ} R^\bot)(X,Y, \xi, \eta)
=
-\widetilde g([J(\overline\nabla _Z A)_\xi, A_\eta]X, Y)
+
\widetilde g([J A_\xi, (\overline\nabla_Z A)_\eta] X, Y).
$$
Hence,  because of (2.3) and (2.5), and using (2.14), we obtain:
$$
(\overline\nabla_{JZ} R^\bot)(X,Y, \xi, \eta)
=
-\widetilde g([(\overline\nabla _Z A)_{J\xi}, A_\eta]X, Y)
+
\widetilde g([A_{J\xi}, (\overline\nabla_Z A)_\eta] X, Y)
=
$$
$$
=
(\overline\nabla_Z R^\bot)(X,Y, J\xi, \eta)
-
2\widetilde g([(\overline\nabla _Z A)_{J\xi}, A_\eta]X, Y).
$$
Lemma is proved.

\section{Proofs of theorems 1, 2.}

{\bf Proof theorem 1.}

Let for some 1-form $\mu$ on $F^{2m},$ the following condition
holds
$$
(\overline\nabla_X b)(Y,Z) = \mu (X) b(Y,Z) \quad \forall X, Y,
Z\in TF^{2m}. \eqno (3.1)
$$
Then for any vector field $\xi \in T^\bot F^{2m},$ we have:
$$
\widetilde g
(
\left(
\overline\nabla_X b
\right)
(Y,Z), \xi )
=
\widetilde g(\mu (X) b(Y,Z), \xi).
$$
Hence,  using (2.1) and (1.6), we obtain:
$$
\widetilde g
(
\left(
\overline\nabla_X A
\right)_\xi Y, Z)
=
\widetilde g (\mu (X) A_\xi Y, Z)
\quad
\forall X, Y, Z\in TF^{2m},
\quad
\forall \xi \in T^\bot F^{2m}.
$$
Thus, the condition (3.1) is equivalent to the condition
$$
\left(
\overline\nabla_X A
\right)_\xi
=
\mu (X) A_\xi,
\quad
\forall X \in TF^{2m},
\quad
\forall \xi \in T^\bot F^{2m}.
\eqno (3.2)
$$
From (3.2) we obtain the equality:
$$
\left(
\overline\nabla_{JX} A
\right)_\xi
=
\mu(JX) A_\xi
\quad
\forall X \in TF^{2m},
\quad
\forall \xi \in T^\bot F^{2m}.
\eqno (3.3)
$$
On the other hand, from (3.2), because of (2.7), we have:
$$
\left(
\overline\nabla_{JX} A
\right)_\xi
= -J
\left(
\mu(X) A_\xi
\right),
\quad
\forall X \in TF^{2m},
\quad
\forall \xi \in T^\bot F^{2m}.
\eqno (3.4)
$$
From (3.3) and (3.4), we obtain:
$$
\mu(JX) A_\xi
= -J
\left(
\mu(X) A_\xi
\right),
\quad
\forall X \in TF^{2m},
\quad
\forall \xi \in T^\bot F^{2m}.
$$
Hence,  for any $Y\in TF^{2m},$ we have:
$$
\mu(JX) A_\xi Y
= - \mu(X)
J
\left(
 A_\xi Y
\right),
\quad
\forall X \in TF^{2m},
\quad
\forall \xi \in T^\bot F^{2m}.
\eqno (3.5)
$$
Using (3.5), we obtain:
$$
\mu(JX)
\widetilde g( A_\xi Y , A_\xi Y )
= - \mu(X)
\widetilde g(
J
\left(
 A_\xi Y
\right),
 A_\xi Y) = 0,
$$
$$
\forall X, Y \in TF^{2m},
\quad
\forall \xi \in T^\bot F^{2m}.
\eqno (3.6)
$$
Since $b\ne 0$ then there exists nondegenerate vector field
$\xi\in T^\bot F^{2m}$, and from (3.6) we come to the equality:
$$
\mu (X) = 0
\quad
\forall X\in TF^{2m}.
$$
Then 1-form $\mu\equiv 0$ and, therefore,
$$
\left(
\overline\nabla_X A
\right)_\xi
= 0,
\quad
\forall X \in TF^{2m},
\quad
\forall \xi \in T^\bot F^{2m}.
\eqno(3.7)
$$
Hence,  because of (2.14), we obtain the conclusion of the
theorem.

\bigskip

{\bf Proof of theorem 2.}

Form (1.3) we obtain:
$$
\nabla_W R(X,Y,Z,V) =
\widetilde g((\overline\nabla_Wb)(X, V), b(Y, Z))+
\widetilde g(b(X, V), (\overline\nabla_Wb)(Y, Z)) -
$$
$$
-
\widetilde g((\overline\nabla_Wb)(X, Z), b(Y, V)) -
\widetilde g(b(X, Z), (\overline\nabla_Wb)(Y, V))
\quad
\forall X,Y,Z,V,W\in TF^{2m}.
$$
Therefore, because of (3.7), $\nabla R\equiv 0$. Theorem is
proved.

\bigskip
{\large\bf References}
\\ \\
$[1]$ Kobayashi S., Nomizu K. Foundations of differential
geometry. Vol. 2. M.: Nauka. 1981.\\
$[2]$ Chen  B.-Y. Geometry of submanifolds. N.-Y.: M. Dekker. 1973.\\
$[3]$ Gray A. Tubes. M.: Mir. 1993.

\end{document}